\documentclass[12pt]{article}
\usepackage{amssymb}
\newcommand{\R}{ \mathbb{R}}

\newcommand{\seq}[1]{(#1_{n})_{n=1,2,\dots}}

\newtheorem{theorem}{Theorem}[section]
\newtheorem{defn}{Definition}[section]
\newtheorem{cor}{Corollary}[section]
\newtheorem{lemma}{Lemma}[section]
\newtheorem{prop}{Proposition}[section]
\newtheorem{rem}{Remark}[section]

\begin{document}

\title{The DFR property for counting processes stopped at an independent random time} 

\author {F. G. Bad\'{\i}a, C. Sang\"{u}esa
 \\ {\small
Departamento de M\'etodos Estad\'{\i}sticos and IUMA. Universidad de
Zaragoza. }
 \\ {\small Zaragoza. Spain. } \\ {\small e-mail: gbadia@unizar.es,csangues@unizar.es} }

\date{}
\begin{titlepage}

\setcounter{page}{1} \maketitle

\bigskip \bigskip
\begin{abstract}
In the present paper we consider general counting processes stopped at a random time $T$, independent of the process.  Provided that $T$ has the decreasing failure rate (DFR) property, we give sufficient conditions on the arrival times so that the number of events occurring
before $T$ preserves the DFR property of $T$.  In particular, when the interarrival times are independent, we consider applications concerning the DFR property of the stationary number of customers waiting in queue for specific queuing models.
\end{abstract}



\end{titlepage}

\section{Introduction} 

Decreasing failure rate (DFR) is a property observed in systems which improve with age.  Although counterintuitive, this property appears in real life.  One of the first examples was observed by Proschan \cite{prtheo} in the study of failures of air-conditioning systems in airplanes. The theoretical explanation of these observed decreasing failure rates is that the data can be described by a mixture model of exponential distributions. Although exponential distributions have a constant failure rate,  their mixtures are DFR.  In fact, the preservation of the DFR property under mixing is well-known.  In addition, mixtures of increasing failure rate distributions may show decreasing failure rates, as shown in Gurland and Sethuraman \cite{guseho} or Finkelstein and Esaulova \cite{fieswh}.

The aim of the present paper is to show that counting processes stopped at a DFR random time inherit this property under rather general assumptions on the interarrival times. As far as
we know, one of the first publications studying the preservation of properties in counting processes observed at a random time (classical renewal processes) is
Esary, Marshall and Proschan \cite{esmash}.  The special case of the Poisson process has
been the subject of much attention: cf. Grandell \cite[Ch. 7]{grmixe},
and the references therein. Ross,
Shanthikumar, and Zhou \cite[Theorem 3.1]{roshon} deal with the preservation of the increasing failure rate property (IFR) in general renewal processes. This subject was also studied in Bad\'{\i}a \cite{german}, together with the preservation of the decreasing reversed hazard rate property.

Several examples of practical applications of counting processes stopped at a random time are described in \cite[p.\ 42]{corene}. For instance, $T$ may be the lifetime of a whole piece of equipment, and we may be concerned with the number
of component replacements, $N(T)$, that are necessary during $T$.
In several queuing models, the stationary length of the queue can be modeled by means of a random variable $N(T)$ (see \cite{hanear} and \cite{assamp}), with $T$ being the stationary waiting time and $N$ being the equilibrium renewal processes of the arrival of customers in the queue. Recent applications have been found in health sciences (see \cite{soesti}), in which the counting process describes the number of tests for a disease (HIV for instance) in a person at risk until the random time $T$ at which this person is infected. Several more applications are described in \cite{basapr} (see also the references therein).

The starting point of the present paper is \cite{BS}, in which it was shown that in a classical renewal process (independent and identically distributed interarrival times) the decreasing failure rate (as well as the logconvexity) is preserved. In Theorem \ref{teasso} we extend this property to possibly dependent interarrival times.  As a particular case, we consider independent and stochastically decreasing interarrival times (Corollary \ref{coinde}) and apply it to queuing models (Section 4).  More specifically, we deduce the DFR property of the stationary number of customers waiting in queue for specific queuing models  in Corollaries \ref{queue} and \ref{c3}. Finally, in the last section we show two examples of counting process with dependent interarrival times in Corollaries \ref{c4} and \ref{c5} which satisfy conditions in Theorem \ref{teasso}.

The next section, Section 2, introduces the main concepts which we are going to use throughout the paper.

\section{Preliminaries}
We consider a general counting process $ \{ N(t): t \geq 0 \} $ defined as follows.  The arrival times in the counting process $ \{ N(t): t \geq 0 \} $ are denoted by $S_n$, $ n=0,1, \ldots $ where $S_0=0$, and the interarrival times are denoted by $\seq{X}$. The arrival and interarrival times satisfy the identity
 \[
 S_n = \sum_{i=1}^n X_i, \quad n=1, 2, \ldots .
 \]
 The only assumption concerning the $X_n$ is that they are nonnegative random variables and are not degenerate at 0. The general counting process is defined through the renewal epochs by the following expression:
 \[
 N(t) = \max \{n: S_n \leq t \},  \quad t \geq 0 .
 \]

Our paper deals with the preservation of the decreasing failure rate property.  First of all, we will recall below the definition of this reliability class (Barlow and Proschan
\cite{baprst}), along with some other definitions of reliability classes that are weaker than the DFR property, which will be used throughout the paper.
 The terms \textquoteleft increasing\textquoteright \hskip
3pt and \textquoteleft decreasing\textquoteright \hskip 3pt mean, as
usual, \textquoteleft non-decreasing\textquoteright \hskip 3pt and
\textquoteleft non-increasing\textquoteright,  \hskip 3pt
respectively.

 \begin{defn} Let $X$ be a nonnegative random variable with $G$ and \mbox{$\overline
 G:=1-G$} the corresponding distribution and reliability
 functions.  $X$  (or $G$) is said to be:
\begin{enumerate} \item \emph{Decreasing failure  rate}  if

$\overline G (z+t)/\overline G(t)$ is increasing in $t$, for all $z\geq 0$;

\item \emph{New worse than used}
(NWU)  if

$\overline{G}(z)\overline{G}(t)\leq\overline{G}(z+t)$, for all
$z,t\geq 0$;

\item \emph{New worse than used in expectation}
(NWUE)  if

\[\mu:=E[X]\leq \int_{0}^{\infty}\overline G (z+t)/\overline G(t)dz,\quad \hbox{for all } t\geq 0;\]

\item \emph{Increasing mean residual life} (IMRL) if
\[m(t):=\int_{0}^{\infty}\overline G (z+t)/\overline G(t)dz,\quad \hbox{ is increasing in } t\geq 0.\]

\end{enumerate}\label{dfgrel}\end{defn}

Note that if the monotonicity and sense of the inequalities in Definition \ref{dfgrel} are reversed, we obtain the dual concepts of \emph{increasing failure rate} (IFR), \emph{new better than used}, \emph{new better than used in expectation} and \emph{decreasing mean residual life}.

\begin{rem} We list some properties of the previous definitions which will be used throughout the paper. \begin{enumerate}
\item
A DFR distribution can show a jump only at the origin, and the interval of support of DFR distributions is $[0,\infty)$ (\cite[p.\ 117]{maolli})).
  \item If $\overline{G}$ is a DFR survival function, and we take $t\geq 0$ fixed, then the function defined by $\bar{H}(z)=\overline{G}(z+t)/\overline{G}(t),\quad z\geq 0$ (residual life survival function) is also a DFR survival function (\cite[p.\ 118]{maolli}).
  \item DFR implies NWU, NWUE, and IMRL \cite[p.\ 181]{maolli}.
\end{enumerate}\label{reages}
\end{rem}

To describe the decreasing failure rate property for counting processes we will use the definition of this property for discrete random variables (see Esary {\it et al.}
\cite{esmash} or Grandell \cite[Ch. 7]{grmixe}), which is recalled in the following.

\begin{defn} Let $X$ be a (possibly defective) nonnegative integer-valued random
variable and $p_{n}$ its corresponding probability mass function
\[p_{n}:=P(X=n)\quad
n=0,1,\dots .
\]
The distribution function is said to be
\emph{discrete decreasing failure  rate} (d-DFR) if $P(X=n)/P(X\geq
n)$ is decreasing in $n$ or, equivalently, if
\begin{eqnarray} P(X\geq n+1)^{2}&\leq&
P(X\geq n)P(X\geq n+2) \quad n=0,1,\dots .
\end{eqnarray}
\label{dfdrel}\end{defn}

Some conditions in our results require the comparison of random variables with respect to the usual stochastic order, which we recall:

\begin{defn} Let $X$ and $Y$ be two random variables with distribution functions $F_{X}$ and $F_{Y}$, respectively. $X$ is said to be \emph{smaller
than $Y$ in the usual stochastic order} (denoted $X\leq_{ST} Y$) if $F_{X}(t)\geq F_{Y}(t)$ for all real  $t$ or, equivalently, if $\overline{F}_{X}(t)\leq\overline{F}_{Y}(t)$ for all real  $t$.\label{dfstor}
\end{defn}

We will write $X=_{ST}Y$ to indicate that the random variables $X$ and $Y$ have the same distribution function.

\begin{rem} The following property of the usual stochastic order
can be found in \cite[p.\ 6]{mustco}, for instance.
Let $X$ and $Y$ be two random variables. $ X \leq_{ST} Y$ if and only if $E[h(X)] \leq ( \geq ) E[h(Y)]$ for all increasing (decreasing) function $h$ such that the expectations exist.

\label{reorde}\end{rem}

   The next concept of association between random variables, introduced by Esary, Proschan and Walkup \cite{esprwa}, will also play an important role in our results. From now on, given two vectors $\mathbf{x}:=(x_{1},\dots x_{n})\in \R^{n}$ and $\mathbf{y}:=(y_{1},\dots ,y_{n})\in \R^{n}$, the notation $\mathbf{x}\leq \mathbf{y}$ will be used to mean $x_{i}\leq y_{i},\ i=1,2,\dots, n$.  A function $f
:\R^{n}\rightarrow \R$ will be said to be \emph{increasing} if  $\mathbf{x}\leq \mathbf{y}$ implies $f(\mathbf{x})\leq f(\mathbf{y})$ (which is the same thing as to say that $f$ is componentwise increasing).

\begin{defn} A random vector $\mathbf{X}:=(X_{1},\dots,X_{n})$ is said to be \emph{associated} if
\begin{equation}Cov(f(\mathbf{X}),g(\mathbf{X}))=E[f(\mathbf{X})\cdotp g(\mathbf{X})]-E[f(\mathbf{X})]\cdotp E[g(\mathbf{X})]\geq0\label{poscor}\end{equation}
for all increasing functions $f,g:\R^{n}\rightarrow \R$ such that the expectation exists.
\end{defn}
\begin{rem} It will be useful throughout the paper to point out that for an associated vector, (\ref{poscor}) also holds true whenever $f$ and $g$ are decreasing. One only has to apply the definition to the increasing functions $-f$ and $-g$. \label{reass1}\end{rem}
\begin{rem} The following interesting properties of association
can be found in \cite[p.\ 681]{maolli} and \cite[p.\ 123]{mustco}, for instance.
\begin{enumerate}

\item If $\mathbf{X}$ is associated and  $f_{i}:\R^{n}\rightarrow\R,\quad i=1,2,\dots,m,\quad m\in \mathbb{N}$ are all increasing (or all decreasing), then $(f_{1}(\mathbf{X}),\dots,f_{m}(\mathbf{X}))$ is associated.
    \item If $\mathbf{X}:=(X_{1},\dots,X_{n})$ has independent components, then it is associated.
    \item The association property is closed under weak convergence.  That is, if we have a sequence of associated random vectors $(\mathbf{X}_{k})_{k=1,2,\dots}$ converging in distribution to $\mathbf{X}$, then $\mathbf{X}$ is associated.
    \end{enumerate}
 \label{reasso}   \end{rem}

Finally, we introduce some notations, concerning conditional distributions, which play an important role in the proofs of our results.  Given a non negative random vector $(\mathbf{X},\mathbf{Y})$, in which $\mathbf{X}$ has dimension $n$ and $\mathbf{Y}$ has dimension $m$,  we will denote by $\{F_{\mathbf{Y}|\mathbf{x}}, \mathbf{x}\in {\mathbb{R}}_{+}^{n}\}$ to a (regular) conditional distribution of $\mathbf{Y}$ given $\mathbf{X}$, that is, a family of distribution functions satisfying for every Borel sets $A\subseteq \R^{n}_{+},\ B\subseteq\R^{m}_{+}$,
\begin{equation}\int_{\{\mathbf{x}\in A,\mathbf{y}\in B\}}dF_{(\mathbf{X},\mathbf{Y})}(\mathbf{x},\mathbf{y})=\int_{\mathbf{x}\in A}dF_{\mathbf{X}}(\mathbf{x})\int_{\mathbf{y}\in B}dF_{\mathbf{Y}|\mathbf{x}}(\mathbf{y}).\label{promed}\end{equation}
 The existence of such family of conditional distribution functions is always guaranteed (cf. \cite[p.107]{kafoun}). Moreover, if in a counting process we consider $\mathbf{Y}=(X_{n+1},\dots, X_{n+m})$, and take $\mathbf{X}$ a random vector of dimension $n$ containing the previous inter arrival times, for each $\mathbf{x}\in {\mathbb{R}}_{+}^{n}$, we will denote by $(Z_{n+1}^{\mathbf{x}},\dots, Z_{n+m}^{\mathbf{x}})$ to a random vector whose distribution function is $F_{(X_{n+1},\dots, X_{n+m})|\mathbf{x}}$.  The family of random vectors $\{(Z_{n+1}^{\mathbf{x}},\dots, Z_{n+m}^{\mathbf{x}}), \  \mathbf{x}\in {\mathbb{R}}_{+}^{n}\}$ will be called a distributional version of $(X_{n+1},\dots, X_{n+m})$ given $\mathbf{X}$. In this way, for every function  $f:\mathbb{R}_{+}^{n}\times \mathbb{R}_{+}^{m}\rightarrow \mathbb{R}_{+}$, the expected value of $f(\mathbf{X},\mathbf{Y})$ can be calculated, thanks to (\ref{promed}), as
   \begin{equation}E[f(\mathbf{X},\mathbf{Y})]=\int_{\mathbb{R}_{+}^{n}}Ef(\mathbf{x},Z_{n+1}^{\mathbf{x}},\dots, Z_{n+m}^{\mathbf{x}})dF_{\mathbf{X}}(\mathbf{x})
   .\label{fubini}\end{equation}

   \section{DFR property under association and decreasing conditions of the interarrival times}\label{second}
In this Section we present our main result concerning the DFR preservation property in a counting process, making use of the conditional distributions of the inter renewal epochs.
  Let $  \{ N(t): t \geq 0 \} $ be a general counting process with inter renewal epochs $\seq{X}$. Let $T$ be a DFR random time independent from the process, whose survival function is denoted by $\overline {F}_T$. Our aim is to show that if $T$ is DFR, then $N(T)$ is d-DFR.  For technical reasons, we will assume that for each $n=1,2,\dots$, the distribution functions of $S_{n}$, the $n$-th renewal epoch, and $T$ don't have any common points of discontinuity. Then the following equality is verified (cf. \cite[p.\ 642]{esmash})
\begin{equation}
P(N(T) \geq n ) = E (\overline {F}_T (S_n)), \quad n =1,2,\dots. \label{nocomm}
\end{equation}

Making use of the previous equality, the d-DFR property for $N(T)$ can be proved using the following lemma, which will be used extensively throughout the paper.

  \begin{lemma}
Let $  \{ N(t): t \geq 0 \} $ be a general counting process with interarrival times $\seq{X}$. Let $T$ be a random time independent of the process. Assume that

\begin{description}
\item {a) } $T$ is a DFR random variable,
\item {b) } $P(T=0)\cdotp P(X_{1}=0)=0$, that is, $T$ and $X_{1}$ don't have simultaneously positive mass at 0.
\end{description}
Then, $N(T)$ is d-DFR if
\begin{eqnarray}E^{2} (\overline {F}_T (X_{1}))&\leq& E [\overline {F}_T (S_{2})];\label{dfrn0}\\
E ^{2}(\overline {F}_T (S_{n+1}))&\leq &E (\overline {F}_T (S_n))E (\overline {F}_T (S_{n+2})),\quad n=1,2,\dots .\label{dfrn}\end{eqnarray}
\label{ledfsu}\end{lemma}
{\bf Proof:} Under the previous assumptions, we have to prove that $N(T)$ is d-DFR, that is,
\begin{equation}
P^{2}(N(T) \geq n +1) \leq P(N(T) \geq n )P(N(T) \geq n+2 ),\quad n=0,1,2\dots . \label{rintin}
\end{equation}
Note firstly that as $X_{1}\leq_{ST} S_{n}$, assumption b) implies that \begin{equation}P(T=0)\cdotp P(S_{n}=0)\leq P(T=0)\cdotp P(X_{1}=0)=0 \quad n=1,2,\dots .\end{equation}
This means that for each $n=1,2\dots $, $S_{n}$ and  $T$ don't have the origin as a common discontinuity point in their distribution functions.  Recalling property 1 in Remark \ref{reages} and that $T$ is DFR, this implies that $S_{n}$ and  $T$  don't have any common discontinuity points in their distribution functions, and hence (\ref{nocomm}) holds.  Thus, taking into account that $P(N(T) \geq 0)=1$, (\ref{rintin}) for $n=0$  becomes
\[ E^{2} (\overline {F}_T (X_{1}))\leq E (\overline {F}_T (S_{2})),\]
whereas, for $n=1,2,\dots$, (\ref{nocomm}) implies that we can rewrite (\ref{rintin}) as
\begin{equation} E ^{2}(\overline {F}_T (S_{n+1}))\leq E (\overline {F}_T (S_n))E (\overline {F}_T (S_{n+2})).\label{equiva}\end{equation}
$\quad \Box$

Next Proposition gives conditions under which (\ref{dfrn0}) is verfied.

\begin{prop} Let $\seq{X}$ be an arbitrary sequence of interarrival times and $T$ a random time independent of them.  Assume that $T$ is a DFR random variable. If $(X_{1},X_{2})$ is associated and $X_{2}\leq_{ST}X_{1}$, then

\[E^2[ \overline {F}_T (X_{1})] \leq E[ \overline {F}_T (S_{2})]\overline {F}_T (0) . \] \label{prasso}

\end{prop}
  {\bf Proof:}  First of all, we note that, as  $\overline {F}_T $ is a DFR survival function, the survival function $ \overline {F}_T(\cdotp)/\overline {F}_T(0) $ satisfies the NWU property (recall 2 and 3 of Remark \ref{reages}). Thus,
\[\frac{\overline {F}_T(x_{1}+x_{2})}{ \overline {F}_T(0)}\geq \frac{\overline {F}_T(x_{1})}{ \overline {F}_T(0)}\frac{\overline {F}_T(x_{2})}{ \overline {F}_T(0)}, \quad x_{1},x_{2}\geq 0\]
and therefore
\[\overline {F}_T(x_{1}+x_{2})\geq \frac{\overline {F}_T(x_{1})\overline {F}_T(x_{2})}{ \overline {F}_T(0)},  \quad x_{1},x_{2}\geq 0. \]
We take $x_{1}=X_{1}$ and $x_{2}=X_{2}$ to write
\[\overline {F}_T(X_{1}+X_{2})\geq \frac{\overline {F}_T(X_{1})\overline {F}_T(X_{2})}{ \overline {F}_T(0)}.  \]
We take expectations in the previous inequality to obtain
\begin{equation}E[\overline {F}_T(S_{2})] \geq \frac{E[\overline {F}_T(X_{1})\overline {F}_T(X_{2})]}{\overline {F}_T(0)}.  \label{primpa}\end{equation}
As $(X_{1},X_{2})$ is associated and the functions $g_{i}:\R^{2}\rightarrow \R, \ i=1,2$ defined by $g_{i}(x_{1},x_{2})=\overline {F}_T(x_{i}),\ i=1,2$ are decreasing, recalling Remark \ref{reass1} we can write
\begin{equation}E[\overline {F}_T(X_{1})\overline {F}_T(X_{2})]\geq E[\overline {F}_T(X_{1})]E[\overline {F}_T(X_{2})]. \label{secpa}\end{equation}
Since $X_{2}\leq_{ST}X_{1}$, we have by Remark \ref{reorde} that
\begin{equation}E[\overline {F}_T(X_{1})]E[\overline {F}_T(X_{2})]\geq E^{2}[\overline {F}_T(X_{1})]. \label{secp2}\end{equation}
The conclusion follows by (\ref{primpa})--(\ref{secp2}). $\quad \Box$

 \medskip In the next result we give sufficient conditions in order to verify (\ref{dfrn}) when we can find distributional versions of $(X_{n+1},X_{n+2})$ given $(X_{1},\dots,X_{n})$ satisfying (\ref{dfrn0}).
\begin{prop}
Let $  \{ N(t): t \geq 0 \} $ be a general counting process with inter renewal epochs $\seq{X}$. Let $T$ be a DFR random time independent from the process. Let $n$ be a fixed natural number and assume that there exists a distributional version of of $(X_{n+1},X_{n+2})$ given $(X_{1},\dots,X_{n})$ (say $\{(Z_{n+1}^{\mathbf{x}},Z_{n+2}^{\mathbf{x}}),\quad  \mathbf{x}\in \R_{+}^{n}\}$) satisfying, for every DFR survival function  $\overline {H}$ with  $\overline {H}(0)=1$, that
\begin{equation}E^{2}[ \overline {H} (Z_{n+1}^{\mathbf{x}})] \leq E[ \overline {H}( Z_{n+1}^\mathbf{x}+Z_{n+2}^{\mathbf{x}})],\quad \mathbf{x}\in N_{n}, \label{concon}\end{equation}
in which $ N_{n}\subseteq R_{+}^{n}$ is such that $P( (X_{1},\dots,X_{n})\in N_{n})=1$. Then,
\begin{equation} E ^{2}(\overline {F}_T (S_{n+1}))\leq E (\overline {F}_T (S_n))E (\overline {F}_T (S_{n+2})). \label{dfrn*}\end{equation}
\label{prdfr2}
\end{prop}

{\bf Proof:}  Note that, without lost of generality, we can assume that (\ref{concon}) is verified for all $ \mathbf{x}\in \mathbb{R}_{+}^{n}$. Outside  $ N_{n}$, we can always define  $(Z_{n+1}^{\mathbf{x}},Z_{n+2}^{\mathbf{x}})=(0,0)$, which obviously satisfies (\ref{concon}).  Thus, in order to check inequality (\ref{dfrn*}), consider $\mathbf{X}:=(X_{1},\dots,X_{n})$, and for each $\mathbf{x}:=(x_{1},\dots,x_{n})\in \mathbb{R}_{+}^{n}$, call $s_{n}:=x_{1}+\dots+x_{n}$.  As $\overline {F}_T $ is a DFR survival function, then the survival function $\overline {H}$ defined as
\begin{equation}\overline {H}(z)=\frac{\overline {F}_T (s_{n}+z)}{\overline {F}_T (s_{n})},\quad z\geq 0\label{auxh}\end{equation}
verifies the DFR property (recall 2 of Remark \ref{reages}). Using (\ref{fubini}), we can write
\begin{eqnarray}E[\overline {F}_T(S_{n+1})]&=&\int_{\mathbb{R}_{+}^{n}}E[\overline {F}_T(s_{n}+Z_{n+1}^{\mathbf{x}})]dF_{\mathbf{X}}(\mathbf{x})\nonumber\\&=&\int_{\mathbb{R}_{+}^{n}}\overline {F}_T (s_{n})E[\overline {H}(Z_{n+1}^{\mathbf{x}})]dF_{\mathbf{X}}(\mathbf{x})\label{dfr1}\end{eqnarray}
As $\overline {H}(0)=1$, we use (\ref{concon}) and Cauchy-Schwartz's inequality to obtain
\begin{eqnarray}&&\int_{\mathbb{R}_{+}^{n}}\overline {F}_T (s_{n})E[\overline {H}(Z_{n+1}^{\mathbf{x}})]dF_{\mathbf{X}}(\mathbf{x})\leq \int_{\mathbb{R}_{+}^{n}}\overline {F}_T (s_{n})E^{\frac{1}{2}}[\overline {H}(Z_{n+1}^{\mathbf{x}}+Z_{n+2}^{\mathbf{x}})]dF_{\mathbf{X}}(\mathbf{x})\nonumber\\&&=\int_{\mathbb{R}_{+}^{n}}\overline {F}_T ^{\frac{1}{2}}(s_{n})E^{\frac{1}{2}}[\overline {F}_T(s_{n}+Z_{n+1}^{\mathbf{x}}+Z_{n+2}^{\mathbf{x}})]dF_{\mathbf{X}}(\mathbf{x})\nonumber\\& &\leq E^{\frac{1}{2}}[\overline {F}_T(S_{n})] E^{\frac{1}{2}}[\overline {F}_T(S_{n+2})].\label{dfr2}\end{eqnarray}
Thus, (\ref{dfr1}) and (\ref{dfr2}) show (\ref{dfrn*}). $\quad \Box$

 \medskip Now we are in a position to prove the main result of this section.
 \begin{theorem}
Let $  \{ N(t): t \geq 0 \} $ be a general counting process with inter renewal epochs $\seq{X}$. Let $T$ be a random time independent from the process. Assume that
\begin{description}
\item {a) }$T$ is a DFR random variable
\item {b) }$T$ and $X_1$ don't have 0 as a common discontinuity point.
\item {c) }$(X_{1},X_{2})$ is associated and $X_{2}\leq_{ST}X_{1}$.
\item {d) }  For each $n=1,2,\dots$ there exists a distributional version of $(X_{n+1},X_{n+2})$ given $(X_{1},\dots,X_{n})$ (say $\{(Z_{n+1}^{\mathbf{x}},Z_{n+2}^{\mathbf{x}}),\quad  \mathbf{x}\in R_{+}^{n}\}$ satisfying, on a set $ N_{n}\subseteq R_{+}^{n}$ such that $P( (X_{1},\dots,X_{n})\in N_{n})=1$, that
    \[(Z_{n+1}^{\mathbf{x}},Z_{n+2}^{\mathbf{x}})\hbox{ is associated and } Z_{n+2}^{\mathbf{x}}\leq_{ST}Z_{n+1}^{\mathbf{x}},\quad \hbox{for all } \mathbf{x}\in N_{n}.\]
    \end{description}
    Then, $N(T)$ is d-DFR. \label{teasso}
    \end{theorem}

{\bf Proof:}   The result will follow by Lemma \ref{ledfsu}.  First of all, (\ref{dfrn0}) follows by condition c) and Proposition \ref{prasso}.
To prove (\ref{dfrn}), fix $n$ and consider $\bar{H}$ a DFR survival function with $\bar{H}(0)=1$.  Condition d), together with Proposition \ref{prasso} applied to $(X_{1},X_{2}):=(Z_{n+1}^{\mathbf{x}},Z_{n+2}^{\mathbf{x}}),\ \mathbf{x}\in N_{n}$ implies that
    \[E^2[ \overline {H} (Z_{n+1}^{\mathbf{x}})] \leq E[ \overline {H} (Z_{n+1}^{\mathbf{x}}+Z_{n+2}^{\mathbf{x}})],\quad \mathbf{x}\in N_{n} . \]
This, together with Proposition \ref{prdfr2}, proves (\ref{dfrn}), and therefore the d-DFR property for $N(T)$.$\quad \Box$
\medskip

As an immediate consequence of Theorem \ref{teasso} we have the following DFR preservation property for independent interarrival times.
\begin{cor}\label{indepe}
Let $\{ N(t): t \geq 0 \}$ be a general counting process with interarrival times $\seq{X}$. Let $T$ be a nonnegative random variable, independent of the process.  Suppose that the following assumptions hold.
\begin{description}
\item {a) }$T$ is a DFR random variable.
\item {b) }$T$ and $X_{1}$ don't have simultaneously positive mass at 0.
\item {c) }The interarrival times are independent.
\item {d) }$X_{n+1}\leq_{ST} X_n$, $n=1,2, \ldots$.
\end{description}
Then, $N(T)$ is d-
DFR.\label{coinde}
\end{cor}

{\bf Proof:}  The result follows as a consequence of Theorem \ref{teasso}.  First of all, due to the independence assumption and using 2 in Remark \ref{reasso}, we have that $(X_{n+1},X_{n+2})$ is associated for each $n=0,1,\dots$.  Therefore, condition c) in  Theorem \ref{teasso} is verified taking $n=0$ and the fact that $X_{2}\leq_{ST} X_1$.  Secondly, due to the independence assumption, a distributional version of $(X_{n+1},X_{n+2})$ given $(X_{1},\dots,X_{n})$ is trivially provided by $\{(Z_{n+1}^{\mathbf{x}}=X_{n+1},Z_{n+2}^{\mathbf{x}}=X_{n+2})\quad  \mathbf{x}\in R_{+}^{n}\}$.  The association property for this vector, in addition to $X_{n+1}\leq_{ST} X_n$, imply condition d) in Theorem \ref{teasso}. $\quad \Box$

\section{Applications in queuing models}
In this section we present some applications of Corollary \ref{coinde} in queuing models.  The first result gives conditions for our preservation property in a pure birth process.
 \begin{cor}
\label {c4}
Let $\{ N(t): t \geq 0 \}$ be a pure birth process with rates $\lambda_n$ increasing in $n$. If $T$ is a DFR random time independent of the process, then $N(T)$ is d-DFR. \label{purbir}
\end{cor}

{\bf Proof:}  A pure birth process is a counting process with independent interarrival times such that the $n$th interarrival time $X_n$ is a exponential random variable with mean $\frac {1}{\lambda_n}$, $n = 1, 2, \ldots$. Then the result is derived from Corollary \ref{indepe}, taking into account that the random variables $X_n$ are decreasing in the usual stochastic order if the sequence $\lambda_n$ is increasing in $n$. $\quad \Box$

\begin{rem} Corollary \ref{purbir} implies, in particular, that a Yule process stopped at a DFR random time is d-DFR.  Recall that a Yule process has rates $\lambda_n=n \lambda$, for some fixed $\lambda>0$. \end{rem}

Also as a consequence of Corollary \ref{coinde}, we give in the next results some applications concerning  preservation results for queuing models. This example makes use of the distributional Little's law \cite{hanear} and uses previous preservation results due to Shanthikumar \cite{shdfrp}. To this end, consider a queuing model which is a FIFO GI/GI/1.  This means that times between the arrivals of customers to the system are i.i.d. random variables,  $\seq{A}$, and the service times $\seq{B}$ are i.i.d. random variables, both being independent processes. Also, there is only one server, which follows the FIFO discipline (first in, first out). We will assume, in addition, that $E[B_{1}]<E[A_{1}]$ (stability of the queue). It is well known that the steady-state number of customers waiting in a queue $L^{*}$ (the total number of customers waiting to be served in a steady-state queue) follows the distributional Little's law (see \cite{hanear} for the original result, or \cite[p.\ 217]{mustco}, for instance, for the specific application we are going to deal with);  that is, $L^{*}=_{ST} N(T)$, in which $\{N(t),\ t\geq 0\}$ is an equilibrium renewal process and $T$ is a random time independent of it.  These quantities are defined in terms of the queue as follows.
\begin{itemize}
\item $\{N(t),\ t\geq 0\}$ is a time stationary version of the renewal counting process describing the arrivals of customers $\seq{A}$, that is, with independent interarrival times. The first interarrival time $X_{1}$ has the equilibrium distribution of $A_{1}$, given by
  \begin{equation}P(X_{1}> t)=\frac{1}{E[A_{1}]}\int_{t}^{\infty}P(A_{1}> u) du,\quad t\geq 0,\label{equili}\end{equation}
    and the following interarrival times $X_{2}$, $X_{3}, \dots$ have the same distribution as $A_{1}$.
    \item $T$ is distributed as the steady-state wating time (the total time waiting in the queue, in the stationary situation).
    \end{itemize}
On the other hand, we will use some results of Shanthikumar (cf. \cite{shdfrp}) giving sufficient conditions for the DFR property of the steady state waiting time $T$.

    We are now in a position to enunciate the following.
    \begin{cor} Consider a FIFO GI/GI/1 queuing system, as described above.  If the interarrival times of customers $\seq{A}$ are NWUE and the service times $\seq{B}$ are DFR random variables, then the steady-state number of customers waiting in queue $L^{*}$ is d-DFR.\label{queue} \end{cor}
{\bf Proof:} The proof is based on Corollary \ref{indepe} and the above representation of $L^{*}$ in terms of $N(T)$, a renewal process with interarrival times $\seq{X}$ stopped at an independent random time $T$, as described above. First of all, the DFR property for the service times implies that $T$, the steady-state wating time, is DFR (see Shanthikumar, \cite[Thm. 5.3]{shdfrp}), and assumption a) in Corollary \ref{indepe} holds true.  To verify assumption b), note that $X_{1}$ has the equilibrium distribution of the arrival epochs of customers, given in (\ref{equili}). Then, $X_{1}$ is continuous, and therefore $X_{1}$ and $T$ don't have 0 as a common discontinuity point in their distribution functions (despite the fact that $T$ is a compound geometric distribution, thus having strictly positive mass at the origin, see for instance, \cite{shdfrp}). The independence assumption c) for the interarrival times $\seq{X}$ is also verified. As for d), it is also known (\cite[p.\ 174]{maolli}, for instance) that if $A_{1}$ is NWUE, then its equilibrium distribution (which is the distribution of $X_{1}$) verifies $A_{1}\leq_{ST}X_{1}$. As noticed before, the remaining interarrival times $X_{2}$, $X_{3}, \dots$ are distributed as  $A_{1}$ and we have therefore $X_{1}\geq_{ST}X_{2}
    =_{ST}X_{3}\dots$, which proves d).   $\quad \Box$

    \begin{rem} In Asmussen {\it et al.}  \cite{assamp}, the asymptotic properties of the number $L$ of customers in a FIFO GI/GI/1 queuing system were studied (note that in Corollary \ref{queue}, we study $L^{*}$, the number of customers in the queue).  Note that if we apply Little's distributional law to $L$ (cf. \cite{assamp}), $L=_{ST}N(T+B)$, where $\{N(t),\ t\geq 0\}$ and $T$ are as above and $B=_{ST}B_{1}$, the service time, and is independent of  $\{N(t),\ t\geq 0\}$ and $T$.  In fact, $T+B$ represent the stationary sojourn time of a customer (the time from arrival to departure from the system).  Note that the conditions in Corollary \ref{queue} don't imply the DFR property for $L$.  In fact, we can ensure the DFR property of both $T$ and $B$, but not the DFR property of $T+B$, as this property is not preserved under convolution. However, as $L^{*}=(L-1)_{+}$, the set of inequalities in Definition \ref{dfdrel} for $L$  and $n=2,3,\dots$ are verified by means of the d-DFR property of $L^{*}$, but they might fail at $n=0$ or $n=1$.\end{rem}

   A similar result to the previous one can be proved for an M/GI/1 FIFO queue.  Recall that an M/GI/1 queuing system is a GI/GI/1 queuing system in which the arrival times follow an exponential distribution.
    \begin{cor}
\label {c3}
Consider an M/GI/1 FIFO queuing system.  If the service times are IMRL, then $L^{*}$ is d-DFR.
\end{cor}

{\bf Proof:} In an M/GI/1 queuing system, the distribution function of the arrival times is exponential and the exponential distribution is NWUE. Then proceeding as in the proof of Corollary \ref{queue}, we conclude  that $L^{*}=N(T)$ satisfies conditions c)--d) of Corollary \ref{coinde}.  Condition  a), that is, the DFR property for $T$, is a consequence of the IMRL property for the service times (see Shanthikumar \cite[Thm. 5.2.]{shdfrp}).   $\quad \Box$

\begin{rem}Recall that the DFR property is stronger than the IMRL and NWUE properties.  Therefore, Corollary \ref{queue} is in particular verified if both the interarrival times and the service times are DFR random variables.  Also, when we consider in Corollary \ref{c3} an M/GI/1 queuing system, we can weaken the DFR condition of the service time to the IMRL condition.\end{rem}

\section{Examples with dependent interarrival times.}
In this Section we show with some examples that Theorem \ref{teasso} can be applied to obtain the DFR property for counting processes with dependent interarrival times. First of all, next Lemma will be useful in order to check the association condition.
\begin{lemma} Let $(Z_{1},Z_{2})$ be a non negative random vector with independent components. Take two increasing (decreasing) functions $g_{i}:\R^{2}_{+}\rightarrow \R, \ i=1,2$. The random vector $(g_{1}(Z_{1},Z_{2}),g_{2}(Z_{1},Z_{2}))$  is associated.
\label{lefunc}
\end{lemma}
{\bf Proof:}  We have, by the independence assumption, that $(Z_{1},Z_{2})$ is associated. (Recall property 2 in Remark \ref{reasso}). On the other hand, if $g_{i}:\R^{2}_{+}\rightarrow \R $ are increasing (decreasing), we can extend them to functions $g^{*}_{i}:\R^{2}\rightarrow \R $ which are also increasing (decreasing), by putting $g^{*}_{i}(z_{1},z_{2})=g_{i}(z_{1}^{+},z_{2}^{+})$, where for any real number $z$, $z^{+}=max(z,0)$. By property 1 in Remark \ref{reasso}, $(g_{1}^{*}(Z_{1},Z_{2}),g_{2}^{*}(Z_{1},Z_{2}))$ is associated and the conclusion follows as $(Z_{1},Z_{2})$ is nonnegative.$\quad \Box$

\medskip
Making use of previous Lemma, we present two examples of counting processes verifying hypothesis in Theorem \ref{teasso}.
\begin{cor}
\label {c4}
Let $ \{ N(t): t \geq 0 \} $ be a general counting process with interarrival times defined as
\[X_{n}=\prod_{i=1}^{n}Y_{i},\quad n=1,2,\dots , \] where $\seq{Y}$ is a sequence of independent nonnegative random variables such that each $ Y_{i}$ has support on $[0,1]$.  If $T$ is a DFR random variable independent of the process and $T$ and $Y_{1}$ don't have 0 as a common discontinuity point, then
\ $N(T)$ is d-DFR. \label{coprod}

\end{cor}
{\bf Proof:} Conditions a) and b) in Theorem \ref{teasso} are satisfied by assumption. For the remaining conditions, let $n=0,1,2,\dots$.  We see by Lemma \ref{lefunc} that $(c_{n}Y_{n+1},c_{n}Y_{n+1} Y_{n+2})$, $ c_{n}\geq0$ is associated.  In particular, for $n=0$, $(X_{1}=Y_{1},X_{2}=Y_{1} Y_{2})$ is associated. Since $Y_{1}  Y_{2}\leq_{ST} Y_{1}$, by the assumption about the support of the $Y_{i}$, condition c) in Theorem \ref{teasso} is satisfied. Consider now $n=1,2,\dots$ and let $c_{n}(\mathbf{x})=\prod_{i=1}^{n}x_{i}, \ \mathbf{x}\in R_{+}^{n}$.  Clearly, a distributional version of $(X_{n+1},X_{n+2})$ given $(X_{1},\dots,X_{n})$ is provided by \[\{(Z_{n+1}^{\mathbf{x}}=c_{n}(\mathbf{x})Y_{n+1},Z_{n+2}^{\mathbf{x}}=c_{n}(\mathbf{x})Y_{n+1} Y_{n+2}),\   \mathbf{x}\in R_{+}^{n}\}.\]  A similar argument as in the case $n=0$, shows that $(Z_{n+1}^{\mathbf{x}},Z_{n+2}^{\mathbf{x}})$ verifies condition d) in Theorem \ref{teasso}. $\quad \Box$
\begin{cor}
\label {c5}
Let $ \{ N(t): t \geq 0 \} $ be a general counting process and let $\seq{Y}$ be a sequence of independent nonnegative random variables. Define the interarrival times by
\begin{equation} X_{n}=\frac{1}{z+Y_{1}+\dots+Y_{n}},\quad n=1,2,\dots,\qquad z\geq 0. \label{express}\end{equation}  Let $T$ be a DFR random variable independent of the process. Then,
\begin{description}
\item {a) } If $z>0$, $N(T)$ is d-DFR.
\item {b) } If $P(Y_{1}>0)=1$ and $z=0$, then $N(T)$ is d-DFR. \end{description} \label{cosum}
\end{cor}
{\bf Proof:} First of all, condition b) in Theorem \ref{teasso} is satisfied by assumption, as $X_{1}=(z+Y_{1})^{-1}$ has no discontinuity at the origin. Secondly, for each $n=0,1,2,\dots$ we see by Lemma \ref{lefunc} that, under the assumptions in parts a) and b), \begin{equation}\left(\frac{1}{c_{n}+Y_{n+1}},\frac{1}{c_{n}+Y_{n+1}+ Y_{n+2}}\right)\quad \hbox{is associated}\quad \hbox{for} \quad c_{n}>0,\label{assocd}\end{equation}
and obviously
\begin{equation}\frac{1}{c_{n}+Y_{n+1}}\geq_{ST}\frac{1}{c_{n}+Y_{n+1}+ Y_{n+2}}, \quad c_{n}> 0.\label{decrea}\end{equation}
Thus, for part a) we conclude from (\ref{assocd}) that
\begin{equation}\left(X_{1}=\frac{1}{z+Y_{1}},X_{2}=\frac{1}{z+Y_{1}+ Y_{2}}\right) \quad \hbox{is associated.}\label{assoc1}\end{equation}
Moreover, $X_{1}\geq_{ST}X_{2}$ by (\ref{decrea}). For part b), $(X_{1}=Y_{1}^{-1},X_{2}=(Y_{1}+ Y_{2})^{-1})$ is also associated using property 3 in Remark \ref{reasso}, as the random vector in (\ref{assoc1}) converges in distribution to $(Y_{1}^{-1},(Y_{1}+ Y_{2})^{-1})$, as $z\downarrow 0$ (in fact, the convergence is almost sure).  As in this case we trivially have $X_{1}\geq_{ST}X_{2}$, condition c) in Theorem \ref{teasso} is satisfied in both cases.
For condition d), observe that, using (\ref{express}), we have
 \[\left(X_{n+1}=\frac{1}{X_{n}^{-1}+Y_{n+1}},X_{n+2}=\frac{1}{X_{n}^{-1}+Y_{n+1}+ Y_{n+2}}\right),\ n=1,2,\dots.\]
We also see by (\ref{express}) that  $P(X_{n}> 0)=1$, in both cases a) and b).  Thus, if we take $N_{n}=\{\mathbf{x}\in R_{+}^{n}|\ x_{n}>0\}$, we have that $P((X_{1},\dots,X_{n}))\in N_{n})=1$. Then, a distributional version of $(X_{n+1},X_{n+2})$ given $(X_{1},\dots,X_{n})$ can be defined by taking
  \begin{equation}\left(Z_{n+1}^{\mathbf{x}}=\frac{1}{x_{n}^{-1}+Y_{n+1}},Z_{n+2}^{\mathbf{x}}=\frac{1}{x_{n}^{-1}+Y_{n+1}+ Y_{n+2}}\right) ,\quad \mathbf{x}\in N_{n}.\label{assoc2}\end{equation}
 We are thus in conditions to apply (\ref{assocd}) and (\ref{decrea}) to conclude that $(Z_{n+1}^{\mathbf{x}},Z_{n+2}^{\mathbf{x}}),\ \mathbf{x}\in N_{n}$ verifies condition d) in Theorem \ref{teasso}. $\quad \Box$

\section*{Aknowledgements}

This work has been supported by
research projects MTM2010-15311 and by FEDER
funds. The first and second authors acknowledge the support of DGA S11 and E64,
respectively.

\end{document}